\documentclass[11pt]{article}

\usepackage[cmtip,arrow]{xy}
\usepackage{amsmath,amssymb,enumerate,pb-diagram,pb-xy}%,pdfsync}
\def\today{\number\day .\number\month .\number\year}
\usepackage[applemac]{inputenc}

\def \1{{\bf 1}}

\def \al{\alpha}

\def \C{{\mathbb C}}
\def \CA{{\cal A}}

\def \CC{{\cal C}}
\def \CD{{\cal D}}
\def \CE{{\cal E}}
\def \CF{{\cal F}}

\def \CS{{\cal S}}

\def \coim{\operatorname{coim}}
\def \coker{\operatorname{coker}}

\def \dto{{\to \hspace{-9pt}\to}}
\def \df{\ \begin{array}{c} _{\rm def}\\ ^{\displaystyle =}\end
{array}\ }

\def \F{{\mathbb F}}

\def \ga{\gamma}

\def \Hom{{\rm Hom}}
\def \Id{{\rm Id}}
\def \im{\operatorname{im}}

\def \ot{\leftarrow}
\def \ph{\varphi}

\def \prf{{\bf Proof: }}
\def \Q{{\mathbb Q}}
\def \qed{\ifmmode\eqno \square 
		\else\noproof\vskip 12pt plus 3pt minus 9pt \fi}
\def \noproof{{\unskip\nobreak\hfill\penalty50\hskip2em\hbox{}%
     \nobreak\hfill $\square$\parfillskip=0pt%
     \finalhyphendemerits=0\par}}
\def \opp{{\rm opp}}

\def \sr{\stackrel}

\def \tot{{\rm tot}}

\def \Z{{\mathbb Z}}
\def \({\left(}
\def \){\right)}
\def \={{\ =\ }}

\renewcommand{\sp}[1]{\left\langle #1\right\rangle}
\newcommand{\ol}[1]{\overline{#1}}

\newtheorem{theorem}{Theorem}[section]

\newtheorem{lemma}[theorem]{Lemma}

\newtheorem{proposition}[theorem]{Proposition}

\begin{document}

\pagestyle{myheadings} \markright{HARMONIC ANALYSIS over ADELIC 
SPACES}

\title{Harmonic analysis over adelic spaces}
\author{Anton Deitmar\footnote{This paper was written during a stay of the author at the Hausdorff Research Institute for Mathematics, Bonn. The author expresses his gratitude for the warm hospitality and the possibility to enjoy the inspiring atmosphere of this place.}}
\date{}
\maketitle

{\bf Abstract:} In \cite{OP}, D. Osipov and A. Parshin developed an approach to harmonic analysis on higher dimensional local fields through categories of filtered vector spaces as in \cite{Osi}.
In the present paper we give a variant of this approach that behaves nicely for inductive arguments.
We establish Pontryagin duality, the Fourier inversion formula, Plancherel formula, and a Poisson summation formula for all dimensions.

$$ $$

\newpage
\tableofcontents

\section*{Introduction}

In the paper \cite{P76} A.N. Parshin developed a construction of adeles for 2-dimensional schemes which later was generalized to arbitrary dimension by by A.A. Beilinson \cite{Beil}.
See \cite{Huber} for a more detailed exposition.
In this calculus, local fields are replaced by higher dimensional local fields which have valuations in higher dimensional groups.
It emerges the necessity to develop a theory of Harmonic Analysis for higher dimensional local fields in order to be able to extend the methods of Tate's thesis \cite{TT} to higher dimensional schemes, as advocated in
\cite{P99}.
In the inspiring paper \cite{OP}, the authors D.V. Osipov and A.N. Parshin present a higher dimensional harmonic analysis based on a sequence of categories of filtered objects $\CC_n$ defined by the first author in \cite{Osi}.
By a careful analysis they proceed from the zero dimensional case, which corresponds to harmonic analysis on finite abelian goups, to the one dimensional case, which corresponds to totally disconnected abelian groups, to the 2-dimensional case.

In the present paper we give a slightly different construction of categories of filtered objects, called $\CA_n$ for distinction.
The approach of this paper preserves nice properties of these categories, for instance, they turn out to be exact categories.
Following \cite{OP}, we define the spaces of smooth functions $\CE$, uniformly smooth functions $\tilde\CE$ and their duals, the spaces of compactly supported distributions $\CE'$ and of compactly supported uniform distributions $\tilde\CE'$.
All these spaces are defined inductively, where the induction runs by the dimension.
Inductively, we define various Fourier transforms on these space, which satisfy the usual inversion laws and Plancherel formulae.

The Poisson summation formula requires a discrete subobject with compact quotient.
If this is given, the outset gives rise to a space $\CD$ of ``smooth functions of compact support'' on which
 there  is a Fourier transform $f\to \hat f$ mapping
$\CD(A)$ to $\CD(\hat A)$, where $\hat A$ is the Pontryaging dual.
We prove  the  Poisson summation formula for functions in $\CD$.

\section{Filtrations}
A partially ordered set $(I,\le)$ can be considered as a category with 
exactly one arrow from $x$ to $y$ if $x\le y$ and no arrow 
otherwise.
In this way partially ordered sets are just the same as small 
categories with $|\Hom(x,y)|\le 1$.
Functors between such categories are the same thing as order 
preserving maps.

A partially ordered set $(I,\le)$ is called \emph{filtering}, if any 
two elements possess lower and upper bounds in $I$, i.e., if for 
$a,b\in I$ there are $x,y\in I$ such that $x\le a,b\le y$.
Viewing $I$ as a category, we also speak of a \emph{filtering category}.

Let $\CA$ be an abelian category.
A \emph{filtration} on an object $A$ of $\CA$ is a functor $F:I_F\to
\CA$ from a filtering category $I_F$ such that
\begin{itemize}
\item $F(\ph)$ is injective (mono) for every arrow $\ph$ in $I_F$,
\item the injective limit of the diagram $F$ is $A$,
\item the projective limit of $F$ is zero.
\end{itemize}
Modulo natural isomorphy this is the same as saying that for each $i
\in I_F$ one gives a sub-object $F(i)$ of $A$ such that $i\le j\ 
\Rightarrow\ F(i)\subset F(j)$ and $\bigcap_iF(i)=0$ as well as $
\bigcup_iF(i)=A$.

If $F$ is a filtration on $A$ and $\ph:B\to A$ a morphism in $\CA$, 
then one can pull back the filtration to get a filtration $\ph^*F$ 
on $B$ by insisting that for each $i$ the diagram
$$
\begin{diagram}
\node{\ph^*F(i)}\arrow{e}\arrow{s}
	\node{F(i)}\arrow{s}\\
\node{B}\arrow{e,t}\ph
	\node{A}
\end{diagram}
$$
be Cartesian.
In particular, if $B$ is a sub-object of $A$ one can write $\ph^F(i)
=F(i)\cap B$.

An order preserving map $\phi:I\to J$ between filtering sets is 
called \emph{cofinal} if for every $j\in J$ there are $i_1,i_2\in I$ 
with
$\phi(i_1)\le j\le\phi(i_2)$.
Let $F:I_F\to\CA$ be a filtration.
A \emph{sub-filtration} is a pair $(\phi,S)$ where $\phi:I_S\to I_F$ 
is a cofinal map and $S:I_S\to\CA$ is a filtration such that the 
diagram of functors
$$
\begin{diagram}
\node{I_S}\arrow{e,t}{\phi}\arrow{se,b}{S}
	\node{I_F}\arrow{s,r}F\\
\node[2]{\CA}
\end{diagram}
$$
commutes up to isomorphy of functors.
In this case we write $F\succcurlyeq S$.
We consider the equivalence relation $\sim$ on the class of all 
filtrations on $\CA$ which is generated by $F\sim S$ whenever $S$ is 
a sub-filtration of $F$.
Then one has $F\sim G$ if and only if there exist filtrations 
$F=F_1,\dots,F_m=G$ such that for each index $j$, either $F_j$ is a 
sub-filtration of $F_{j+1}$ or the other way round.

\begin{lemma}
For two filtrations of an object $A$ the following are equivalent.
\begin{enumerate}[\rm (a)]
\item $F\sim G$.
\item There is a filtration $H$ such that $H\succcurlyeq F$ and $H
\succcurlyeq G$.
\item For every $i\in I_F$ there are $j_1,j_2\in I_G$ such that
$$
G(j_1)\ \subset\ F(i)\ \subset\ G(j_2),
$$
and the same with reversed roles of $F$ and $G$.
\end{enumerate}
\end{lemma}

\prf (a) $\Rightarrow$ (b).
This follows, if we can show that if $F$ and $G$ have a common 
sub-filtration $S$, then there exists $H$ with $H\succcurlyeq F,G$.
For this let $I_H$ be the disjoint union of $I_F$ and $I_G$.
Define $H$ on objects by $F$ or $G$ whichever is appropriate.
Next define a partial order on $I_H$ by 
$$
i\le j\quad\Leftrightarrow\quad H(i)\subset H(j),
$$
where, properly speaking, the inclusion means the existence of an 
injection which commutes with the chosen injections to $A$.
These chosen injections then make up the images of morphisms under 
$H$.
The existence of a common sub-filtration implies that $F$ and $G$ are 
indeed sub-filtrations of $H$.
The converse direction (b) $\Rightarrow$ (a) is trivial.
Finally, (c) is a reformulation of (b).
\qed

We will write an equivalence class of filtrations as $(A,[F])$, 
where $A$ is the direct limit of $F$, which does not depend on the 
choice of the representative $F$.

\begin{proposition}
Every filtration $F$ with countable index set $I_F$ is equivalent to 
a filtration with index set $\Z$.
\end{proposition}

\prf
This is clear as a countable filtering set $I$ admits a cofinal map 
$\Z\to I$.
\qed

From now on we will restrict to countable filtrations only.

\section{The strong category}
\subsection{The category $\CA_n$}
Let $\CA_0$ be a full abelian subcategory of $\CA$ which is closed 
under isomorphy, i.e., if $c\in\CA_0$ and $a\in\CA$ is isomorphic to 
$c$, then $a\in\CA_0$.
We construct a sequence $\CA_0,\CA_1,\dots$ of categories of 
filtered objects in $\CA$ as follows.
Firstly, we view each object $A$ of $\CA_0$ as trivially filtered 
with $I_F$ consisting of two elements $0$ and $\infty$ and $F(0)=0$ 
as well as $F(\infty)=A$.

For the induction assume $\CA_{n-1}$ already defined as a category 
of certain classes of filtered objects in $\CA$ and the morphisms 
are certain morphisms in $\CA$.
We define the objects of $\CA_{n}$ to be the equivalence
 classes of countable filtrations $(A,[F])$ in $\CA$ together
  with a class of filtrations $[E_{i,j}]$ on the quotient 
$F(j)/F(i)$  for each $i\le j\in I_F$ such that $(F(j)/F(i),[E_
{i,j}])$ is an object of $\CA_n$ and such that the natural maps
$$
F(i)/F(j)\ \to\ F(k)/F(l)
$$
are morphisms in $\CA_{n-1}$ whenever $k\ge i$ and $l\ge j$.
So, strictly speaking, an object of $\CA_n$ is an object of $\CA$ 
with a filtration and with filtrations on all quotients and again 
filtrations on all of their quotients and so forth.
We will not write out all the filtrations, they will be implicit in 
saying that an object belongs to $\CA_n$. 

Let $(A,[F])$ and $(B,[G])$ be objects of $\CA_{n}$.
A morphism in $\CA_{n}$ is a  morphism $\phi:A\to B$ such that for 
every $i\in I_F$ and  every $j\in I_G$ 
there exist $i_0\le i$ and $j_0\ge j$ such that for every $i'\le i_0
$ and every $j'\ge j_0$ one has 
$\phi(F(i'))\subset G(j)$ and $\phi(F(i))\subset G(j')$, and the 
induced map
$$
F(i)/F(i')\ \to\ G(j')/G(j)
$$
is a morphism in $\CA_{n-1}$.

\begin{theorem}\label{exactCat}
For every $n$, the category $\CA_n$ is an additive category which 
contains finite limits.

The category $\CA_n$ is in general not abelian.
\end{theorem}

\prf
This is clear for $n=0$.
For $n\ge 1$ we start by showing additivity.
Let $\phi,\psi$ be morphisms in $\CA_n$ from $(A,[F])$ to $(B,G])$.
In $\CA$, we can form the sum $\phi +\psi: A\to B$.
We have to show that this gives a morphism in $\CA_n$.
Let $i\in I_F$ and $j\in I_G$.
There are $i_\phi,i_\psi\le i$ and $j_\phi,j_\psi\ge j$ such that
for every $i'\le i_\phi,i_\psi$ and every $j\ge j_\phi,j_\psi$ the 
morphisms induced by $\phi$ and $\psi$,
$$
F(i)/F(i')\ \to\ G(j')/G(j)
$$
are morphisms in $\CE_{n-1}'$.
Choose $i_0\le i_\phi,i_\psi$ and $j_0\ge j_\phi,j_\psi$, then for 
every $i'\ge i_0$ and every $j'\ge j_0$ the morphism induced by $
\phi+\psi$ from $F(i)/F(i')$ to $G(j')/G(j)$ is in $\CA_{n-1}$.
This implies that $\phi+\psi$ is a morphism of $\CA_n$, so this 
category is closed under addition of morphisms.
Next for products.
Let $(A,[F])$ and $(B,[G])$ be in $\CA_n$.
The product $A\times B$ exists in $\CA$.
Define a filtration $F\times G$ by $I_{F\times G}=I_F\times I_G$ 
with the product order, i.e., $(i,i')\le (j,j')$ is equivalent to $i
\le j$ and $i'\le j'$. 
Then $I_{F\times G}$ is filtering.
Define $F\times G(i,j)=F(i)\times G(j)$.
This gives a filtration on $A\times B$.
We give $(A\times B,[F\times G])$ a canonical structure of an object 
of $\CA_n$ as follows.
For $(i,i')\le (j,j')$ we have an $\CA$-isomorphism 
\begin{eqnarray*}
F\times G(j,j')/F\times G(i,i')&=&(F(j)\times G(j'))/(F(i)\times G
(i'))\\
&\cong& F(j)/F(i)\times G(j')/G(i').
\end{eqnarray*}
The product filtration on the right hand side will now make $(A
\times B,[F\times G])$ an object of $\CA_n$.
Next we show that it is indeed a product.
The projections $p_A,p_B: A\times B\to A,B$ are in $\CA_n$.
The universal property follows from the one in $\CA$.
As $\CA$ is abelian, $A\times B$ also has the coproduct property in 
$\CA$.
It is straightforward to see that the same holds in $\CA_n$.
So $\CA_n$ is an additive category.

Since we have products and coproducts, the existence of finite 
limits will follow from the existence of kernels and cokernels.
For kernels let $\phi:(A,[F])\to (B,[G])$ be a morphism.
Let $\al:K\to A$ be the kernel of $\phi$ in $\CA$.
Equip $K$ with the filtration $H$ induced from the embedding $K
\hookrightarrow A$, so $H(i)=F(i)\cap K$.
Then $H(j)/H(i)=F(j)\cap K/F(i)\cap K$ injects into $F(j)/F(i)$.
We equip $H(j)/H(i)$ with the filtration induced by this injection 
and so forth.
In this way $(K,[H])$ becomes an object of $\CA_n$ and the embedding 
$K\hookrightarrow A$ is a morphism in $\CA_n$.
Let $\beta:(Z,[J])\to (A,[F])$ be a morphism in $\\CE_n'$ with $\phi
\circ\beta =0$.
We have the diagram
$$
\begin{diagram}
\node{(Z,[J])}\arrow{s,l,..}{\gamma}\arrow{se,t}\beta
\arrow{see,t}0\\
\node{(K,[H])}\arrow{e,t}\al
	\node{(A,[F])}\arrow{e,t}\phi
		\node{(B,[G]).}
\end{diagram}
$$
As $K$ is the kernel of $\phi$ in $\CA$, there exists a $\ga:Z\to K$ 
making the diagram commute.
We have to show that $\ga$ is in $\CA_n$.
So let $i\in I_J$ and $j\in I_H=I_F$. Then there are $i_\beta\le i$ 
and $j_\beta\ge J$ such that
$$
\begin{array}{ccccc}
\beta(J(i'))&\subset& F(j)&&\forall i'\le i_\beta\\
\beta(J(i))&\subset& F(j')&&\forall j'\ge j_\beta
\end{array}
$$
and the morphisms
$$
J(i)/J(i')\ \to\  F(J')/F(j)
$$
are in $\CA_{n-1}$.
Now $\beta$ factorizes over $\ga$, so $J(i)/J(i')$
 maps into the sub-object $F(j')\cap K/F(j)\cap K$ and 
as the filtration on $F(j')\cap K/F(j)\cap K$ is induced from $F
(j')/F(j)$, the map $J(i)/J(i')\to F(j')\cap K/F(j)\cap K$ is in $
\CA_{n-1}$.
It follows that $\ga$ is in $\CA_n$, 
i.e. the category $\CA_n$ possesses kernels.

The existence of cokernels follows by reversing all arrows.
We only give the definition of the filtration on a cokernel.
Let $\phi:(A,[F])\to (B,[G])$ be a morphism in $\CA_n$ and let $
\delta:B\to C$ be a cokernel in $\CA$.
The filtration $H$ on $C$ is defined as $I_H=I_G$ and $H(i)=\delta(G
(i))$. 

It remains to give an example in which $\CA_n$ is not abelian.
Take a field $k$ and let $\CA=\CF=\CA_0$ be the category of all $k$-
vector spaces.
We show that $\CA_1$ is not abelian by giving a morphism with 
trivial kernel and cokernel which is not an isomorphism.
Let $V\in\CA$ of infinite $k$-dimension.
Let $F$ be the filtration of all finite dimensional subspaces and 
let $G$ be the filtration of all subspaces.
Then the identity map $(V,[F])\to (V,[G])$ is in $\CA_1$, has 
trivial kernel and cokernel, but, as $F$ and $G$ are not equivalent, 
it is not an isomorphism in $\CA_n$. 
\qed

A sequence in $\CA_n$,
$$
\begin{diagram}
\node 0\arrow{e}
	\node{A}\arrow{e,t}{\al}
		\node{B}\arrow{e,t}{\beta}
			\node{C}\arrow{e}
				\node 0
\end{diagram}
$$
is called \emph{exact}, if $\al$ is the kernel of $\beta$ and $\beta
$ is the cokernel of $\al$.
By a \emph{kernel} we mean a map which is the kernel of its 
cokernel.
Likewise, a \emph{cokernel} is the cokernel of its kernel.

\begin{proposition}
$\CA_n$ with the class of sequences which are exact in $\CA_n$, is 
an exact category.
\end{proposition}

\prf
The only non-trivial point is to show that the pullback of a 
cokernel is a cokernel.
So let $\phi:A\to B$ be a cokernel in $\CA_n$ and let $\ph:C\to B$ 
be an arbitrary map in $\CA_n$.
We have a Cartesian diagram
$$
\begin{diagram}
\node{(P,[J])}\arrow{e,t}{\phi'}\arrow{s}
	\node{(C,[H])}\arrow{s,r}{\ph}\\
\node{(A,[F])}\arrow{e,t}{\phi}
	\node{(B,[G])}
\end{diagram}
$$
We want to show that $\phi'$ is a cokernel.
As $\CA$ is an abelian category, $\phi'$ is surjective, so we only 
need to show that the filtration on $C$ is induced by $\phi'$.
For this recall that $P$ is the kernel of the map $(\phi-\ph):A
\times C\to B$, so the filtration on $P$ is induced from the product 
filtration on $A\times C$ and $\phi'$ is derived from the projection 
$A\times C\to C$.
Let $i\in I_H$.
Then there is $j\in I_G$ such that $\ph(H(i))\subset G(j)$.
As $\phi$ is a cokernel, we have $G(j)=\phi(F(j))$.
Now we claim that $H(i)=\phi'((F(j)\times H(i))\cap P)$.
Trivially the right hand side is contained in $H(i)$.
For the converse direction we can assume that $\CA$ is a subcategory 
of the category of modules of a ring, which means that we can use 
elements.
So let $x\in H(i)$.
As $\phi(H(i))\subset G(j)=\phi(F(j))$, there exists $y\in F(j)$ 
such that $\phi(y)=\ph(x)$, so $(y,x)\in (F(j)\times H(i))\cap P$ 
and this proves the claim.
\qed

Let $\phi:(A,[F])\to (B,[G])$ be a kernel.
We say that the filtration $F$ is \emph{induced} by $G$, if $I_F=I_G$ and for every $i\in I_F$ one has $F(i)=\phi^{-1}(I_G)$.

If $\phi:(A,[F])\to (B,[G])$ is a cokernel, then we say that $G$ is \emph{induced from $F$}, if $I_G=I_F$ and for every $i\in I_G$ one has $G(i)=\phi(F(i))$.

\begin{lemma}
If $\phi:(A,[F])\to (B,[G])$ is a (co-)kernel, then up to $\CA_n$-isomorphism one can choose the filtration $F$ $(G)$ as induced from $G$ $(F)$ and one can assume that for $i\le j$ the induced map $F(j)/F(i)\to G(j)/G(i)$ is a (co-)kernel in $\CA_{n-1}$.
One can assume that the filtrations on $F(j)/F(i)$ and $G(j)/G(i)$ are induced one by the other accordingly, and so forth.
\end{lemma}

\prf
A (co-)kernel as constructed in the proof of Theorem \ref{exactCat} has this property.
Any two (co-)kernels of the same map are isomorphic.
\qed

\subsection{The category $\CS_n$}
A morphism $\phi$ in $\CA_n$ is called \emph{quasi-strong}, if the 
canonical map
$$
\coim(\phi)\ \to\ \im(\phi)
$$ 
is an isomorphism.
Here, as usual,
$$
\coim(\phi)\=\coker(\ker(\phi)),\quad\mbox{and}\quad \im(\phi)\=\ker
(\coker(\phi)).
$$
A morphism is quasi-strong if and only if it can be written as a cokernel followed by a kernel.
This factorization is unique up to isomorphism.

\begin{lemma}
Isomorphisms are quasi-strong and the composition of two quasi-strong maps is 
quasi-strong.
\end{lemma}

\prf
Isomorphisms are clearly quasi-strong.
We prove the second assertion by induction.
It is clear for $n=0$.
For $n>0$ let $\phi:A\to B$ and $\psi:B\to C$ be quasi-strong.
Let $K_\phi$ and $K_\psi$ be their kernels, then $\phi$ factorizes 
into a cokernel followed by a kernel as
$$
A\ \dto\ A/K_\phi\ \hookrightarrow\ B
$$
and likewise for $\psi$.
Write $X$ for $A/K_\phi$.
We have a diagram
$$
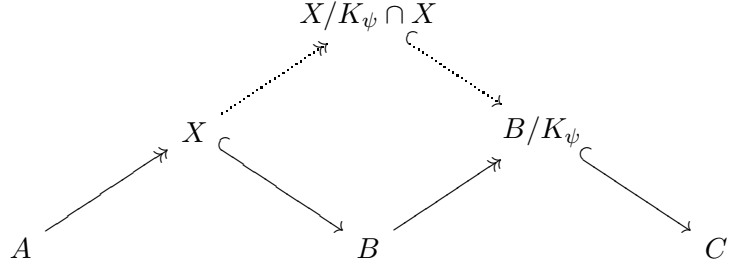
\begin{diagram}
\node[3]{X/K_\psi\cap X}\arrow{se,..,J}\\
\node[2]{X}\arrow{ne,..,A}\arrow{se,J}
	\node[2]{B/K_\psi}\arrow{se,J}\\
\node{A}\arrow{ne,A}
	\node[2]{B}\arrow{ne,A}
		\node[2]{C}
\end{diagram}
$$
The dotted arrows exist in $\CA$.
On $X/K_\psi\cap X$ we have two filtrations, one induced by the 
embedding into $B/K_\psi$ and one induced by the projection from $X
$.
As the filtrations on $X$ and $B/K_\psi$ can both be assumed to be 
induced by one filtration on $B$, it turns out that the two 
filtrations on $X/K_\psi\cap X$ can be assumed to agree.
Taking quotients of the various filtrations, one sees that the 
middle square of the diagram iterates, and so one can deduce that 
indeed the dotted arrows are a kernel and a cokernel respectively in 
$\CA_n$.
The claim follows.
\qed

We define the \emph{strong category} $\CS_n$ as a subcategory of $\CA_n$ inductively as follows.
For $n=0$ we set $\CS_0=\CA_0$.
For $n>0$ we call
an object $(A,[F])$ a \emph{strong object}, if for any $i\le j$ 
the quotient $F(j)/F(i)$ is strong in $\CA_{n-1}$ and for any $i
\le j\le k$ the sequence
$$
F(j)/F(i)\to F(k)/F(i)\to F(k)/F(j)
$$ 
is exact in $\CA_{n-1}$.
A morphism $\phi:A\to B$ is \emph{strong}, if $A$ and $B$  are strong, and $\phi$ is quasi-strong.
The \emph{strong category} $\CS_n$ is the category of strong objects and strong morphisms in $\CA_n$.

\subsection{Completion}
For an object $(A,[F])$ of $\CS_n$ we define the completion $\bar A$ 
in $\CS_n$ together with an injection $A\hookrightarrow 
\bar A$ in $\CS_n$ inductively.
The map $A\to\bar A$ is an endofunctor of $\CS_n$, which is a 
projection in the sense that the given injection $\bar A\to \bar{\bar 
A}$ is an isomorphism.

For $n=0$ we define $\bar A=A$ and the injection is the identity 
map.
For $n>0$ we define
$$
\bar A\=
\lim_{^\to_j}\lim_{^\ot_i}\ol{F(j)/F(i)}.
$$
The filtration $\bar F$ on $\bar A$ is defined by
$$
\bar F(j)\=\lim_{^\ot_i}\ol{F(j)/F(i)}.
$$
To see that this defines an object of $\CS_n$ we have to find a 
natural $\CS_{n-1}$-structure on $\bar F(j)/\bar F(i)$.
We get this by showing that there is a natural isomorphism $\bar F
(j)/\bar F(i)\cong \ol{F(j)/F(i)}$ as part of the next proposition.

\begin{proposition}\label{ExactCompl}
For $j\le k$ we have a natural isomorphism 
$$
\bar F(k)/\bar F(j)\ \cong\ \ol{F(k)/F(j)}.
$$
The completion functor $\CS_n\to\CS_n$ is well-defined and exact.
\end{proposition}

\prf
Note that the assertions are independent of the ambient abelian category $\CA$.
So we can enlarge $\CA$ and assume that it is the full module category of a commutative ring with unit.

All assertions of the proposition are clear if $n=0$.
We will prove these assertions together by an inductive argument.
So assume them proven for $n-1$.
For $(A,[F])$ in $\CS_n$ and $i\le j\le k$ consider the exact sequence in 
$\CS_{n-1}$,
$$
0\to F(j)/F(i)\to F(k)/F(i)\to F(k)/F(j)\to 0.
$$
By induction hypothesis the sequence
$$
0\to \ol{F(j)/F(i)}\to \ol{F(k)/F(i)}\to \ol{F(k)/F(j)}\to 0
$$
is exact for every $i\le j$.
As the last item in the sequence does not depend on $i$, we can take the projective limit over $i$ to get an exact sequence in $\CA$,
$$
0\to \bar F(j)\to \bar F(k)\to \ol{F(k)/F(j)}\to 0.
$$
This gives the first claim and defines the completion 
functor on $\CS_n$.
Let 
$$
\begin{diagram}
\node 0\arrow{e}
	\node{A}\arrow{e,t}{\al}
		\node{B}\arrow{e,t}{\beta}
			\node{C}\arrow{e}
				\node 0
\end{diagram}
$$
be an exact sequence in $\CS_n$.
We have to show that the sequence
$$
\begin{diagram}
\node 0\arrow{e}
	\node{\bar A}\arrow{e,t}{\bar\al}
		\node{\bar B}\arrow{e,t}{\bar\beta}
			\node{\bar C}\arrow{e}
				\node 0
\end{diagram}
$$
is exact in $\CS_n$.
For this recall that the filtrations $F$ and $H$ on $A$ and $C$ are 
induced by the filtration $G$ on $B$ and all filtrations are 
countable, which means that we can assume $I_F=I_G=I_H=\Z$ and $F(i)
=\al^{-1}(G(i))$ as well as $H(i)=\beta(G(i))$.
Therefore we get an exact sequence in $\CS_{n-1}$,
$$
0\to F(j)/F(i)\to G(j)/G(i)\to H(j)/H(i)\to 0
$$
for $i\le j$.
By the induction hypothesis the sequence
$$
0\to \ol{F(j)/F(i)}\to \ol{G(j)/G(i)}\to \ol{H(j)/H(i)}\to 0
$$
is exact.
The functor of taking projective limits is left exact, so we get an 
exact sequence in $\CA$,
$$
0\to \lim_{^\ot_i}\ol{F(j)/F(i)}\to \lim_{^\ot_i}\ol{G(j)/G(i)}\to 
\lim_{^\ot_i}\ol{H(j)/H(i)}\to R^1\lim_{^\ot_i}\ol{F(j)/F(i)}.
$$
The last item denotes the first right derived functor of $\lim_{\ot}
$.
We claim that the last map is also surjective in $\CA$.
For this we need a lemma.

\begin{lemma}
For $i'\le i$ the natural map $\ol{F(j)/F(i')}\to \ol{F(j)/F(i)}$ is 
surjective.
\end{lemma}

\prf
The map $F(j)/F(i')\to F(j)/F(i)$ is surjective in $\CS_{n-1}$, 
hence by induction hypothesis the lemma follows.
\qed

This Lemma implies that the projective system $(\ol{F(j)/F(i)})_i$ 
satisfies the Mittag-Leffler condition.
As $\CA$ is the module category of a ring it follows that
$$
R^1\lim_{^\ot_i}\ol{F(j)/F(i)}\=0.
$$
(See, for example, Proposition 1 in \cite{Roos}.)
From this it follows that the sequence
$$
0\to \lim_{^\ot_i}\ol{F(j)/F(i)}\to \lim_{^\ot_i}\ol{G(j)/G(i)}\to 
\lim_{^\ot_i}\ol{H(j)/H(i)}\to 0
$$
is exact in $\CA$.
Taking direct limits, we see that the sequence $0\to\bar A\to\bar B
\to\bar C\to 0$ is exact in $\CA$.
As the filtrations on both sides are the induced ones, it is also 
exact in $\CS_n$.
\qed

The injection $A\to\bar A$ comes by taking limits of the maps $F(j)/
F(i)\to\ol{F(j)/F(i)}$.
A morphism $\al:A\to B$ in $\CS_n$ naturally induces a morphism on 
the completions $\bar\al:\bar A\to\bar B$.
An object $A$ of $\CS_n$ is called \emph{complete}, if the natural 
map $A\to\bar A$ is an isomorphism.

\section{Pontryagin dual}
We will now specialize to $\CA$ being a module category of a ring.
So let $R$ be a commutative ring with unit and let $\CA$ be the 
category of $R$-modules.
Let $\CA_0$ be the subcategory of finite modules, i.e., those, which 
are finite as sets.
We define a functor $\hat\cdot : \CS_n^\opp\to\CS_n$ together with a 
natural transformation $\delta:\Id\to\hat{\hat\cdot}$ as follows.
For $n=0$ let
$$
\hat A\=\Hom_\Z(A,\Q/\Z).
$$
This is the Pontryagin dual.
Then $\hat A$ is an $R$-module through the rule $r\al(a)=\al(ra)$ 
for $a\in A$ and $r\in R$.
Further the map $\delta:A\to\hat{\hat A}$ given by $\delta(a)(\al)=
\al(a)$ is an isomorphism by the Theorem of Pontryagin.

Next suppose that $\hat\cdot$ is already defined for $\CS_{n-1}$.
For an object $(A,[F])$ of $\CS_n$ we define
$$
\hat A\=\lim_{^\to_i}\lim_{^\ot_j}\widehat{F(j)/F(i)}.
$$
Then $\hat A$ has a filtration $\hat F$ with $I_{\hat F}=I_F^\opp$ 
the same set with opposite order and
$$
\hat F(i)\=\lim_{^\ot_j}\widehat{F(j)/F(i)}.
$$
As in Proposition \ref{ExactCompl} one sees that $\hat F(i)/\hat F
(j)=\widehat{F(j)/F(i)}$ and hence $\hat\cdot$ is a well defined 
functor.
By definition one gets $\hat{\hat{A}}\cong\bar A$ and the map $
\delta$ is the natural injection.
So in particular, if $A$ is complete, then $\delta$ is a natural 
isomorphism $A\to\hat{\hat A}$.

\begin{proposition}
The functor $\hat\cdot$ from $\CS_n^\opp$ to $\CS_n$ is exact.
\end{proposition}

\prf Similar to the proof of Proposition \ref{ExactCompl}.
\qed

\subsection{Compact and discrete objects}
We define \emph{compact objects} of $\CS_n$ as follows.
For $n=0$, every object of $\CS_0$ is compact.
For $n>0$ an element $(A,[F])$ is called compact if $A$ is complete, there is $j$ with $F(j)=A$, and every quotient $F(j)/F(i)$ is compact in $\CS_{n-1}$.

Dually we define the notion of \emph{discrete objects}.
Every object of $\CS_0$ is discrete.
For $n>0$, an object $(A,[F])$ is called discrete if there exists $i$ with $F(i)=0$ and every quotient $F(j)/F(i)$ is discrete in $\CS_{n-1}$.

\begin{proposition}
Let $A$ be an object of $\CS_n$.
If $A$ is compact, then $\hat A$ is discrete.
If $A$ is discrete, then $\hat A$ is compact.
\end{proposition}

\prf
An easy induction.
\qed

\begin{proposition}
Let $K,D$ be subobjects of $A\in\CS_n$, where $D$ is discrete and $K$ is compact.
Then $D\cap K$ is finite.
\end{proposition}

\prf
This is clear for $n=0$.
Let $n>0$.
There is $j$ such that $K\subset F(j)$ and there is $i\le j$ such that $D\cap F(i)=0$.
Therefore $D\cap K$ injects into $F(j)/F(i)$ and is the intersection of a discrete and a compact subobject of $F(j)/F(i)$, hence the claim follows by induction hypothesis.
\qed

\subsection{Subobjects}
In this section we show that subobjects in $\CS_n$ are the same as submodules of the ring $R$.

\begin{lemma}
Let $(A,[F])$ be in $\CS_n$ and let $T\subset A$ be a submodule.
Then, up to isomorphy, there is a unique structure of an $\CS_n$-object on $T$ such that the injection $T\hookrightarrow A$ is a kernel.
\end{lemma}

\prf
Uniqueness is clear, since kernels are uniquely determined up to isomorphy.
The claim is clear for $n=0$.
Let $n>0$.
On $T$ fix the filtration $F_T(j)=T\cap F(j)$.
Then the quotient $T\cap F(j)/T\cap F(i)$ injects into $F(j)/F(i)$, thus has a unique $\CS_{n-1}$ structure making the injection a kernel.
For $i\le j\le k$ one has the commutative diagram
$$
\begin{diagram}
	\node{F_T(j)/F_T(i)}\arrow{e}\arrow{s,J}
		\node{F_T(k)/F_T(i)}\arrow{e}\arrow{s,J}
			\node{F_T(k)/F_T(j)}\arrow{s,J}\\
	\node{F(j)/F(i)}\arrow{e,J}
		\node{F(k)/F(i)}\arrow{e,A}
			\node{F(k)/F(j)}			
\end{diagram}
$$
The lower row is exact, the verticals are kernels, hence the upper row also is exact.
\qed

\section{Smooth functions}
We keep $\CA$ equal to the category of $R$-modules and $\CS_0$ the category of finite modules.
For any $A\in\CA$ let $C(A)$ be the complex vector space of all 
maps from $A$ to $\C$.
For a morphism $\phi: A\to B$ in $\CA$, we get the pullback $\phi^*:C(B)\to C(A)$ defined by $\phi^*(\ph)=\ph\circ\phi$.
If $A$ and $B$ are in $\CS_0$, we also get a push-forward $\phi_*:C(A)\to C(B)$ defined by
$$
\phi_*(\ph)(x)\=\sum_{y:\phi(y)=x}\ph(y),
$$
where the empty sum is interpreted as zero.

Note that if $\phi$ is an injective morphism in $\CA$, then the definition of $\phi_*$ also makes sense and defines $\phi_*:C(A)\to C(B)$.

For $A,B\in\CA_0$ we also define
$$
\phi_!(\ph)(x)\=\frac {|\coker(\phi)|}{|\ker(\phi)|}\sum_{y:\phi(y)=x}\ph(y).
$$

\begin{lemma}
For any two composable morphisms in $\CS_0$ one has $(\psi\phi)_*=\psi_*\phi_*$ and $(\psi\phi)_!=\psi_!\phi_!$.
\end{lemma}

\prf
To fix notations, suppose $\phi:A\to B$ and $\psi:B\to C$.
Then
\begin{eqnarray*}
\psi_*\phi_*f(x)&=& \sum_{b:\psi(b)=x}\sum_{a:\phi(a)=b} f(a)\\
&=& \sum_{a:\psi(\phi(a))=x}f(a)\=(\psi\phi)_*f(x).
\end{eqnarray*}
In the case of lower shriek, one gets the same identity with the factors $\frac {|\coker(\phi)||\coker(\psi)|}{|\ker(\phi)||\ker(\psi)|}$ and $\frac {|\coker(\phi\psi)|}{|\ker(\phi\psi)|}$ respectively, so that $(\phi\psi)_!=\phi_!\psi_!$ is equivalent to $$
|\coker(\psi\phi)||\ker(\phi)||\ker(\psi)|=|\ker(\phi\psi)||\coker(\phi)||\coker(\psi)|.
$$
We denote by $K_\phi, C_\phi, I_\phi$ the kernel, cokernel and image of $\phi$ and likewise for $\psi$.
We get exact sequences
$$
0\to K_\phi\to A\to I_\phi\to 0
$$
and 
$$
0\to I_\phi\to B\to C_\phi\to 0.
$$
Analogous sequences holds for $\psi$ and $\psi\phi$, giving the following identities
\begin{eqnarray*}
|A|&=& |K_\phi||I_\phi|\= |K_{\psi\phi}||I_{\psi\phi}|\\
|B|&=& |K_\psi||I_\psi|\= |I_\phi||C_\phi|\\
|C|&=& |I_\psi||C_\psi|\= |I_{\psi\phi}||C_{\psi\phi}|.
\end{eqnarray*}
These imply the claim by an easy computation.
\qed

Following \cite{OP}, for each $n\ge 0$ we now define two functors  $\CE_n$ and $\tilde\CE_n$ from $\CS_n^\opp$ to the category of complex vector spaces as follows.
For $A\in\CS_0$ we define $\CE_0(A)=\tilde\CE_0(A)=C(A)$ and $\CE_0(\phi)=\tilde\CE_0(\phi)=\phi^*$ as above.
Now suppose $\CE_{n-1}$ and $\tilde\CE_{n-1}$ already defined, then for an object $(A,[F])$ of $\CS_n$ we define the ``space of smooth functions'' as
$$
\CE_n(A,[F])\ \df\ \lim_{^\ot_j}\lim_{^\to_i}\CE_{n-1}(F(j)/F(i)),
$$
and the ``space of uniformly smooth functions'' as
$$
\tilde\CE_n(A,[F])\ \df\ \lim_{^\to_i}\lim_{^\ot_j}\tilde\CE_{n-1}(F(j)/F(i)),
$$
were the limits are taken with respect to $\pi_{ijk}^*$ and $\al_{ijk}^*$.
Let $\phi:(A,[F])\to (B,[G])$ be a kernel or cokernel.
Then the filtrations can be assumed one induced by the other and for $i\le j$ the resulting map $\phi_{ij}:F(j)/F(i)\to G(j)/G(i)$ again a kernel or cokernel respectively. 
The maps $\phi^*:\CE_n(B)\to\CE_n(A)$ and $\phi^*:\tilde\CE_n(B)\to\tilde\CE_n(A)$ is then defined as the limits of the $\phi_{ij}$.

Dually, we define functors of ``distributions'' as follows. 
These are functors $\CE_n'$ and $\tilde\CE_n'$ from $\CS_n$ to the category of $\C$-vector spaces.
Again, $\CE_0'(A)=\tilde\CE_0'(A)=C(A)$ and $\CE_0'(\phi)=\tilde\CE_0'(\phi)=\phi_*$.
For $n>0$ we define
$$
\CE_n'(A,[F])\ \df\ \lim_{^\to_j}\lim_{^\ot_i}\CE_{n-1}'(F(j)/F(i)),
$$
and
$$
\tilde\CE_n'(A,[F])\ \df\ \lim_{^\ot_i}\lim_{^\to_j}\tilde\CE_{n-1}'(F(j)/F(i)),
$$
where the limits are taken with respect to $\pi_*$ and $\al_*$.

\begin{lemma}\label{4.1}
Let $\phi: (A,[F])\to (B,[G])$ be a cokernel in $\CS_n$ then the map $\phi^*:\CE_n(B)\to\CE_n(A)$ is injective.
\end{lemma}

\prf
Clear for $n=0$.
For $n>0$ and $i\le j\le k$ we get a commutative diagram by induction hypothesis,
$$
\begin{diagram}
\node{\CE_{n-1}(F(k)/F(i))}
	\node{\CE_{n-1}(G(k)/G(i))}\arrow{w,t,J}{\phi^*}\\
\node{\CE_{n-1}(F(k)/F(j))}\arrow{n,r,J}{\pi^*}
	\node{\CE_{n-1}(G(k)/G(j)).}\arrow{w,t,J}{\phi^*}\arrow{n,r,J}{\pi^*}
\end{diagram}
$$
Taking injective limits with injective connection morphisms preserves injectivity, therefore the induced map
$$
\phi^*:\lim_{^\to_i}\CE_{n-1}(G(j)/G(i))\ \to\ \lim_{^\to_i}\CE_{n-1}(F(j)/F(i))
$$
is injective for every $j$.
Taking projective limits is a left exact functor, so the map $\phi^*:\CE_n(B)\to\CE_n(A)$ is injective.
\qed

\begin{lemma}
There are natural perfect pairings of complex vector spaces
$$
\sp{\cdot,\cdot}: \CE_n'(A,[F])\times \CE_n(A,[F])\ \to\ \C,
$$
and
$$
\sp{\cdot,\cdot}: \tilde\CE_n'(A,[F])\times \tilde\CE_n(A,[F])\ \to\ \C.
$$
For a strong morphism $\phi:A\to B$ in $\CS_n$ these satisfy
$$
\sp{\phi_*f,g}\=\sp{f,\phi^*g},
$$
if $f\in \CE_n'(B)$ or $\tilde\CE_n'(B)$ and $g\in \CE_n(A)$ or $\tilde\CE_n(A)$.
\end{lemma}

\prf
For $n=0$ the pairing on $C(A)\times C(A)$ is given by
$$
\sp{f,g}\=\sum_{a\in A}f(a)g(a).
$$
This sets up an isomorphism $C(A)\cong C(A)^*$.
The induction comes from the fact that the dual space of an injective limit is the projective limit of the duals and vice versa.
\qed

Let $(A,[F])$ be an object of $\CS_n$.
We define a map $(t_a)_*:\CE_n'(A)\to \CE_n'(A)$ inductively such that $
(t_{a+a'})_*=(t_a)_*(t_{a'})_*$ and such that for every strong 
morphism $\phi:A\to B$ the diagram
$$
\begin{diagram}
\node{\CE_n'(A)}\arrow{e,t}{(t_a)_*}\arrow{s,r}{\phi_*}
	\node{\CE_n'(A)}\arrow{s,r}{\phi_*}\\
\node{\CE_n'(B)}\arrow{e,t}{(t_{\phi(a)})_*}
	\node{\CE_n'(B)}
\end{diagram}
$$
commutes.
For $n=0$ one sets $(t_a)_*f(x)=f(x-a)$ and the claim follows from a 
computation.
For $n>0$ let $i\le j\le k$ and assume that $F(k)$ contains $a$.
By induction hypothesis the diagram
$$
\begin{diagram}
\node{\CE_{n-1}'(F(k)/F(i))}\arrow{e,t}{(t_a)_*}
\arrow{s,r}{(\pi_{ijk})_*}
	\node{\CE_{n-1}'(F(k)/F(i))}\arrow{s,r}{(\pi_{ijk})_*}\\
\node{\CE_{n-1}'(F(k)/F(j))}\arrow{e,t}{(t_a)_*}
	\node{\CE_{n-1}'(F(k)/F(j))}
\end{diagram}
$$
commutes.
If $a$ is contained in $F(j)$ then also the diagram
$$
\begin{diagram}
\node{\CE_{n-1}'(F(j)/F(i))}\arrow{e,t}{(t_a)_*}
\arrow{s,r}{(\al_{ijk})_*}
	\node{\CE_{n-1}'(F(j)/F(i))}\arrow{s,r}{(\al_{ijk})_*}\\
\node{\CE_{n-1}'(F(k)/F(i))}\arrow{e,t}{(t_a)_*}
	\node{\CE_{n-1}'(F(k)/F(i))}
\end{diagram}
$$
commutes.
Thus we can take limits to get a map $(t_a)_*: \CE_n'(A)\to \CE_n'(A)$.
The claimed properties of $(t_a)_*$ follow inductively.
The same notion is used for the analogous maps on $\tilde\CE_n'(A)$.

On the other hand we similarly define maps 
$t_a^*:\CE_n(A)\to \CE_n(A)$ and likewise on $\tilde\CE_n(A)$  such that $t_{a+a'}^*=t_a^*t_{a'}^*$ and 
such that for every strong morphism $\phi:A\to B$ the diagram
$$
\begin{diagram}
\node{\CE_n(A)}\arrow{e,t}{t_a^*}
	\node{\CE_n(A)}\\
\node{\CE_n(B)}\arrow{e,t}{t_{\phi(a)}^*}\arrow{n,r}{\phi^*}
	\node{\CE_n(B)}\arrow{n,r}{\phi^*}
\end{diagram}
$$
commutes.

\begin{proposition}
If $A\ne 0$, the space $\CE_n(A)^A$ of $A$-invariants in $\CE_n(A)$ 
is one-dimensional.
Every strong morphism $\phi:A\to B$ in $\CS_n$ induces a non-zero map $\phi^*: \CE_n(B)^B\to \CE_n(A)^A$.
There is a canonical basis element $\1_A\in\CE_n(A)^A$ with $\phi^*(\1_B)=\1_A$ for every strong $\phi$.

The analogous assertions hold for $\tilde\CE_n(A)$.
\end{proposition}

\prf
For $n=0$ the invariants are just the constant functions, which implies the claim.
The canonical element is the constant function of value $1\in\C$.

For $n>0$ one gets
$$
\CE_n(A)^A\=\lim_{^\ot_j}\(\lim_{^\to_i} \CE_{n-1}(F(j)/F(i))\)^{F(j)}\=\lim_{^\ot_j}\lim_{^\to_i} \CE_{n-1}(F(j)/F(i))^{F(j)}.
$$
This is a limit over one dimensional spaces, hence the dimension of $\CE_n(A)$ is at most one.
As all the maps that make up the limits are non-zero, the space is non-zero.
The functoriality follows by induction.
\qed

\begin{lemma}
For $(A,[F])\in \CS_n$ there is  a natural injective linear map $\tau:\CE_n(A)\hookrightarrow C(A)$ such that for every strong morphism $\phi:A\to B$ in $\CS_n$ the diagram
$$
\begin{diagram}
\node{\CE_n(A)}\arrow{e,t,J}{\tau}
	\node{C(A)}\\
\node{\CE_n(B)}\arrow{n,r}{\phi^*}\arrow{e,t,J}{\tau}
	\node{C(B)}\arrow{n,r}{\phi^*}
\end{diagram}
$$
commutes.
Further, $\tau$ commutes with the $A$-translation action, i.e., for every $a\in A$ one has $t_a^*\tau=\tau t_a^*$.
Likewise, there is an analogous map $\tilde\tau:\tilde\CE_n(A)\hookrightarrow C(A)$. 
\end{lemma}

\prf
For $n=0$ the map $\tau$ is the identity map and the assertions are clear.
For $n>0$ and $i\le j\le k$, using Lemma \ref{4.1} one gets commutative diagrams
$$
\begin{diagram}
\node{\CE_{n-1}(F(k)/F(i))}\arrow{e,t,J}{\tau}
	\node{C(F(k)/F(i))}\\
\node{\CE_{n-1}(F(k)/F(j))}\arrow{n,r,J}{\pi_{ijk}^*}\arrow{e,t,J}{\tau}
	\node{C(F(k)/F(j))}\arrow{n,r,J}{\pi_{ijk}^*}
\end{diagram}
$$
and
$$
\begin{diagram}
\node{\CE_{n-1}(F(j)/F(i))}\arrow{e,t,J}{\tau}
	\node{C(F(j)/F(i))}\\
\node{\CE_{n-1}(F(k)/F(i))}\arrow{n,r}{\al_{ijk}^*}\arrow{e,t,J}{\tau}
	\node{C(F(k)/F(i)).}\arrow{n,r}{\al_{ijk}^*}
\end{diagram}
$$
Therefore one can define $\tau:\CE_n(A)\to C(A)$ as the limit of those maps.
Taking injective limits with injective connection maps preserves injectivity and taking projective limits is left exact, therefore $\tau$ is indeed injective.
Now let $\phi: (A,[F])\to (B,[G])$ be a kernel or cokernel.
By induction hypothesis for $i\le j$ the diagram
$$
\begin{diagram}
\node{\CE_{n-1}(F(j)/F(i))}\arrow{e,t,J}{\tau}
	\node{C(F(j)/F(i))}\\
\node{\CE_{n-1}(G(j)/G(i))}\arrow{n,r}{\phi^*}\arrow{e,t,J}{\tau}
	\node{C(G(j)/G(i))}\arrow{n,r}{\phi^*}
\end{diagram}
$$
commutes.
Taking limits the claimed diagram commutes.
The last assertion is clear.
\qed

\section{Fourier transform on $\CE$}

\subsection{Definition of $\F$}
We consider $\CE_n$ and $\tilde\CE_n$ as contravariant functors on $\CS_n$ and their duals $\CE_n'$ and $\tilde\CE_n'$ as covariant functors.
We will define several Fourier transforms, i.e., natural isomorphisms of functors as follows:
\begin{eqnarray*}
\F:\CE_n &\to& \tilde\CE_n'\circ\hat\cdot\\
\F':\CE_n' &\to& \tilde\CE_n\circ\hat\cdot\\
\tilde\F:\tilde\CE_n &\to& \CE_n'\circ\hat\cdot\\
\tilde\F':\tilde\CE_n' &\to& \CE_n\circ\hat\cdot.
\end{eqnarray*}
For $A\in\CS_n$ we have to define a map $\F=\F_A:\CE_n(A)\to \tilde\CE_n'(\hat A)$ with the property that for every strong morphism $\phi:A\to B$ the diagram
$$
\begin{diagram}
\node{\CE_n(A)}\arrow{e,t}{\F}
	\node{\tilde\CE_n'(\hat A)}\\
\node{\CE_n(B)}\arrow{e,t}{\F}\arrow{n,r}{\phi^*}
	\node{\tilde\CE_n'(B)}\arrow{n,r}{\hat\phi_*}
\end{diagram}
$$
commutes.
We start with $n=0$.
We define
$$
\F(f)(\al)\=\frac 1{|A|}\sum_{a\in A}f(a)e^{-2\pi i\al(a)}.
$$
To show the desired property in this case let $\phi:A\to B$ be a morphism in $\CS_0$ and let $f\in \CE_0(B)$.
Then for $\al\in\hat A$,
\begin{eqnarray*}
\F\phi^*f(\al) &=& \frac 1{|A|}\sum_{a\in A} f(\phi(a)) e^{-2\pi i\al(a)}\\
&=& \frac 1{|A|}\sum_{b\in B}f(b)\sum_{a:\phi(a)=b}e^{-2\pi i\al(a)}.
\end{eqnarray*}
If $b$ lies in the image of $\phi$, then choose $a_0$with $\phi(a_0)=b$.
The second sum becomes
$$
e^{-2\pi i\al(a_0)}\sum_{a\in\ker\phi}e^{-2\pi i\al(a)},
$$
which is zero unless $\al\in(\ker\phi)^\perp=\im\hat\phi$.
Therefore, $\F\phi^*f(\al)$ is zero unless $\al=\hat\phi(\beta_0)$ for some $\beta_0$, in which case it equals
\begin{eqnarray*}
\frac{|\ker\phi|}{|A|}\sum_{b\in\im\phi}f(b)e^{-2\pi i\al(a_0)}
&=& \frac{|\ker\phi|}{|\coker\phi||A|}\sum_{b\in B} f(b)e^{-2\pi i\beta_0(b)}\sum_{\beta\in\ker\hat\phi}e^{-2\pi i\beta(b)}\\
&=& \frac 1{|B|}\sum_{b\in B}f(b)\sum_{\beta:\hat\phi(\beta)=\al}e^{-2\pi i\beta(b)}\\
&=& \frac 1{|B|}\sum_{\beta:\hat\phi(\beta)=\al}\F f(\beta)\=\hat\phi_*\F f(\al),
\end{eqnarray*}
which is the desired identity.

For $n>0$ we define $\F$ as follows.
Let $(A,[F])\in\CS_n$.
We assume that $\F$ has already been defined on $\CS_{n-1}$, so for $i\le j\le k$ there are commutative diagrams
$$
\begin{diagram}
\node{\CE_{n-1}(F(k)/F(i))}\arrow{e,t}{\F}
	\node{\tilde\CE_{n-1}'(\hat F(i)/\hat F(k))}\\
\node{\CE_{n-1}(F(k)/F(j))}\arrow{e,t}{\F}\arrow{n,r}{\pi_{ijk}^*}
	\node{\tilde\CE_{n-1}'(\hat F(j)/\hat F(k))}\arrow{n,r}{(\al_{kji})_*}
\end{diagram}
$$
and
$$
\begin{diagram}
\node{\CE_{n-1}(F(j)/F(i))}\arrow{e,t}{\F}
	\node{\tilde\CE_{n-1}'(\hat F(i)/\hat F(j))}\\
\node{\CE_{n-1}(F(k)/F(i))}\arrow{e,t}{\F}\arrow{n,r}{\al_{ijk}^*}
	\node{\tilde\CE_{n-1}'(\hat F(i)/\hat F(k)).}\arrow{n,r}{(\pi_{kji})_*}
\end{diagram}
$$
Note that $\hat\al_{ijk}=\pi_{kji}$ and $\hat\pi_{ijk}=\al_{kji}$.
This allows us to take limits to obtain
$$
\F: \CE_n(A)\ \to\ \tilde\CE_n'(\hat A).
$$
For the functorial property let $\phi:(A,[F])\to (B,[G])$ be a kernel or cokernel and assume that the filtrations $F$ and $G$ are induced one by the other.
For $i\le j$ there is a commutative diagram
$$
\begin{diagram}
\node{\CE_{n-1}(F(j)/F(i))}\arrow{e,t}\F
	\node{\tilde\CE_{n-1}'(\hat F(i)/\hat F(j))}\\
\node{\CE_{n-1}(G(j)/G(i))}\arrow{e,t}{\F}\arrow{n,r}{\phi^*}
	\node{\tilde\CE_{n-1}'(\hat G(i)/\hat G(j)).}\arrow{n,r}{\hat\phi_*}
\end{diagram}
$$
Taking limits we get
$$
\begin{diagram}
\node{\CE_n(A)}\arrow{e,t}{\F}
	\node{\tilde\CE_n'(\hat A)}\\
\node{\CE_n(B)}\arrow{e,t}{\F}\arrow{n,r}{\phi^*}
	\node{\tilde\CE_n'(\hat B)}\arrow{n,r}{\hat\phi_*}
\end{diagram}
$$
as claimed.

The definitions of $\F'$, $\tilde \F$, and $\tilde\F'$ are completely analogous.
We only have to fix the definitions for $n=0$.
So let $A\in\CS_0$ and $f\in C(A)$.
Then
$$
\F'f(\al)\df \sum_{a\in A}f(a) e^{2\pi i\al(a)},
$$
and $\tilde\F f=\F f$ as well as $\tilde\F'f=\F'f$.

\begin{theorem}
(Inversion formula)\\
The transformations $\F$ and $\tilde\F'$ are inverse to each other in the sense that there are canonical isomorphisms $\F\tilde\F'(\tilde\CE_n')\cong \tilde\CE_n'$ and $\tilde\F'\F(\CE_n)\cong\CE_n$.
In the same sense, $\F'$ and $\tilde\F$ are inverses of each other.
\end{theorem}

\prf
The claim holds for $n=0$ and follows in general by induction.
\qed

\subsection{Plancherel formula}

\begin{theorem}
(Plancherel formula)\\
Let $A\in\CS_n$.
For every $f\in\CE_n(A)$ and every $g\in\CE_n'(A)$ one has
$$
\sp{f,g}\=\sp{\F f,\F' g}.
$$
For $f\in\tilde\CE_n(\hat A)$ and $g\in\tilde\CE_n'(\hat A)$ one has
$$
\sp{f,g}\=\sp{\tilde\F f,\tilde\F' g}.
$$
\end{theorem}

\prf
For $n=0$  this is the Plancherel formula for finite abelian groups.
For $n>0$ let 
$$
f\in\CE_n(A)\=\lim_{^\ot_j}\lim_{^\to_i}\CE_{n-1}(F(j)/F(i)),
$$
and 
$$
g\in\CE_n'(A)\=\lim_{^\to_j}\lim_{^\ot_i}\CE_{n-1}'(F(j)/F(i)).
$$
Then
$$
\F f\in\tilde\CE_n'(\hat A)\=\lim_{^\ot_j}\lim_{^\to_i}\tilde\CE_{n-1}'(\hat F(i)/\hat F(j)),
$$
and 
$$
\F' g\in\tilde\CE_n(\hat A)\=\lim_{^\to_j}\lim_{^\ot_i}\tilde\CE_{n-1}(\hat F(i)/\hat F(j)).
$$
Suppose that $i\in I_F$ is sufficiently small and $j\in I_F$ sufficiently large.
Then the components of $f$ and $g$,
$$
f_{ij}\in\CE_{n-1}(F(j)/F(i)),\quad\mbox{and}\quad g_{ij}\in\CE_{n-1}'(F(j)/F(i))
$$ 
both exist and by induction hypothesis we get
$$
\sp{f,g}\=\sp{f_{ij},g_{ij}}\=\sp{\F f_{ij},\F' g_{ij}}.
$$
We further can insist that $i$ is small enough and $j$ big enough such that $\F f_{ij}=(\F f)_{ij}$ and $\F' g_{ij}=(\F' g)_{ij}$ as well as $\sp{(\F f)_{ij},(\F' g)_{ij}}=\sp{\F f,\F' g}$.
This is the first claim.
The second follows in a similar way.
\qed

\subsection{Fourier transform through functions}
In the case $n=1$ any $(A,[F])\in\CS_1$ can be equipped with he topology generated by the sets of the form $a+F(i)$ for $a\in A$ and $i\in I_F$.
Then $A$ is a totally disconnected group which is locally compact if $A$ is complete.
Then we have a Haar measure and can identify the spaces $\CE,\CE'$ with the spaces of locally constant functions on $A$ and its dual space, the space of compactly supported distributions.
The latter space contains the space $\CD(A)$ of locally constant functions of compact support, on which a Fourier transform is defined via the Haar measure.
In this section we will prove (Theorem \ref{functions}), that this Fourier transform coincides with our given one.

\begin{lemma}
Let $A\in\CS_n$.
There is a canonical pairing
$$
(\cdot,\cdot):A\times\hat A\ \to\ \Q/\Z
$$
such that for every strong morphism $\phi$ one has $(\phi(a),\al)=(a,\hat\phi(\al))$.
\end{lemma}

\prf
For $n=0$ define $(a,\al)=\al(a)$ and the claim follows.
For $n>0$ let $a\in A$ and $\al\in\hat A$.
Then there are $i\le j$ such that $a\in F(j)$ and $\al\in\hat F(i)=\lim_{^\ot_k}\widehat{F(k)/F(i)}$.
Let $\al_j$ be the projection of $\al$ to $\widehat{F(j)/F(i)}$ and let $a_i$ be the projection of $a$ to $F(j)/F(i)$, then 
$$
(a,\al)\df (a_i,\al_j)
$$
does not depend on the choice of $i,j$.
To show the claimed property, let $\phi:(A,[F])\to (B,[G])$ be a kernel or cokernel, where we assume that thhe filtrations are induced one by the other through $\phi$.
For $i\le j$ the map $\phi$ induces $\phi_{ij}:F(j)/F(i)\to G(j)/G(i)$.
For suitable indices $i\le j$ one gets
\begin{eqnarray*}
(\phi(a),\beta)&=& (\phi(a)_i,\beta_j)\=(\phi_{ij}(a_i),\beta_j)\\
&=& (a_i,\widehat{\phi_{ij}}\beta_j)\\
&=& (a_i,\hat\phi(\beta)_j)\= (a,\hat\phi(\beta)),
\end{eqnarray*}
as claimed.
\qed

So in particular, every $\al\in\hat A$ defines an element $e^{2\pi i\al(\cdot)}$ of $C(A)$.

\begin{lemma}
For $f\in\CE_n(A)$ and $\al\in\hat A$ there exists a unique element $fe^{2\pi i\al}$ in $\CE_n(A)$ such that 
$$
\tau(fe^{2\pi i\al})(x)\= \tau(f)(x)e^{2\pi i\al(x)}.
$$
\end{lemma}

\prf
The uniqueness is clear by the injectivity of $\tau$.
We prove existence.
The claim is trivial for $n=0$.
For $n>0$ let $f$ be an element of $\CE_n(A)=\lim_{^\ot_j}\lim_{^\to_i}\CE_{n-1}(F(j)/F(i))$.
Let $\al\in\hat A$, say $\al\in\hat F(i)=\lim_{^\ot_j}\widehat{F(j)/F(i)}$, where we are free to decrease $i$ if necessary.
Let $a\in A$, say $x\in F(j)$, and let $f_j$ be the projection of 
$f$ to $\lim_{^\to_i}\CE_{n-1}(F(j)/F(i))$, say $f_j\in\CE_{n-1}(F
(j)/F(i))$ with the same $i$ as above.
Let $x_i$ be the projection of $x$ to $F(j)/F(i)$, then $\tau(f)(x)=\tau(f_j)(x_i)$.
Let $\al_j$ be the projection of $\al$ to $\widehat{F(j)/F(i)}$.
Then there is $f_je^{2\pi i\al_j}\in\CE_{n-1}(F(j)/F(i))$ such that 
$$
\tau(f_je^{2\pi i\al_j})(x_i)\=\tau(f_j)(x_i)e^{2\pi i\al_j(x_i)}\=\tau(f)(x)e^{2\pi i\al(x)}.
$$
Define the element $fe^{2\pi i\al}$ of $\CE_n(A)$ by the components $f_je^{2\pi i\al_j}$.
The claim follows.
\qed

\begin{theorem}\label{functions}
For $f\in\CE_n'(A)$ one has
$$
\tau(\F'f)(\al)\=\sp{\1_A,fe^{2\pi i\al}}.
$$
Likewise, for $g\in\tilde\CE_n'(A)$ one has
$$
\tau(\tilde\F'g)(\al)\=\sp{\1_A,ge^{2\pi i\al}}.
$$
\end{theorem}

\prf
The claims are clear for $n=0$.
For $n>0$ let $f$ be an element of $\CE_n'(A)=\lim_{^\to_j}\lim_{^\ot_i}\CE_{n-1}'(F(j)/F(i))$.
Then $\F'f$ lies in $\tilde\CE_n(\hat A)=\lim_{^\to_j}\lim_{^\ot_i}\tilde\CE_{n-1}(\widehat{F(j)/F(i)})$.
Now for suitable indices, $(\F'f)_{ij}=\F'f_{ij}$ and
$$
\tau(\F'f)(\al)\=\tau(\F'f_{ij})(\al_{ij})\=\sp{\1_{F(j)/F(i)},f_{ij}e^{2\pi i\al_{ij}}}\=\sp{\1_A,fe^{2\pi i\al}}.
$$
\qed

\section{The Poisson summation formula}
For given $(A,[F])$ in $\CS_n$ we can assume that the index set of the filtration $F$ is $\Z$.
The filtration on each quotient $F(i+1)/F(i)$ can also be assumed to be indexed by $\Z$ and so forth.
In this way we get a \emph{total filtration} $F_\tot$ on $A$, indexed by $\Z^n$ with the lexicographic ordering.

The filtration $F$ also gives rise to a filtration $\hat F$ on the Pontryagin dual $\hat A$ where the ordering of the index set is turned around and
$\hat F(i)/\hat F(j)\cong\widehat{F(j)/F(i)}$.
This implies that the filtration on $F(i+1)/F(i)$ induces a filtration on $\hat F(i)/\hat F(i+1)$ and so on.
In this way we get a total filtration $\hat F_\tot$ on $\hat A$.
As confusion is unlikely, we will write $F$ for $F_\tot$.
For instance, for $z=(z_1,\dots,z_n)\in\Z^n$ one has $F(z_1)\subset F(z)$.

\begin{proposition}
For $z\in\Z^n$ one has
$$
\hat F(z)\=\{\chi\in\hat A: \chi(F(z))=0\}\= F(z)^\perp.
$$
\end{proposition}

\prf
The claim is immediate for $n=0$.
For $n>0$ let $\chi\in\hat F(z)$.
Note that $F(z_1)\subset F(z)\subset F(z_1+1)$, so that $\hat F(z_1+1)\subset\hat F(z)\subset \hat F(z_1)$ and hence $\chi\in\hat F(z_1)=\lim_{^\ot_j}\widehat{F(j)/F(z_1)}$.
This implies $\chi(F(z_1))=0$.
Therefore, $\chi$ induces an element $\bar\chi$ of $\widehat{F(z_1+1)/F(z_1)}$.
Let $a\in F(z)\subset F(z_1+1)$ and let $\bar a$ be the projection of $a$ to $F(z_1+1)/F(z_1)$.
Then $\bar a\in \bar F(z_2,\dots,z_n)$, where $\bar F$ is the filtration on $F(z_1+1)/F(z_1)$.
Further $\chi(a)=\bar\chi(\bar a)$ and the latter is zero by induction hypothesis.
It follows $\chi\in F(z)^\perp$, hence we have shown ``$\subset$''.

For the other direction let $\chi$ lies in $F(z)^\perp\subset F(z_1)^\perp$.
This means that $\chi\in\lim_{^\ot_j}\widehat{F(j)/F(z_1)}=\hat F(z_1)$.
Let $\bar\chi$ be the projection of $\chi$ to $\widehat{F(z_1+1)/F(z_1)}$.
Then $\bar\chi$ lies in $F(z_2,\dots,z_n)^\perp$ and the latter equals $\hat F(z_2,\dots,z_n)$ by induction hypothesis.
This implies the claim.
\qed

Let
\begin{equation}\label{S}
0\ \to\ D\ \to\ A\  \to\ K\ \to\ 0
\end{equation}
be an exact sequence in $\CS_n$ such that $D$ is a discrete object and $K$ a compact one.
When we choose a filtration $F$ on $A$ we consider $D$ and $K$ equipped with the induced filtrations which we write $D\cap F$ and $F_K$ respectively.

We want to define a Fourier transform $f\mapsto\hat f$ on a certain space of functions $f$ on $A$.
For a locally compact group one has to choose a Haar measure in order to define a Fourier transform.
In general, there is no canonical choice of Haar measure, except for discrete groups (counting measure), and compact groups (normalized measure).
This means that for a locally compact group $A$ an exact sequence (\ref{S}) with a discrete $D$ and a compact $K$ gives rise to a canonical Haar-measure, and thus a canonical Fourier transform.

For any subset $S\subset A$ we define $\1_S$ to be the indicator function of the set $S$, i.e., it takes the value $1$ on $S$ and $0$ outside $S$.
We define the space $\CD_n(A)\subset C(A)$ to be the linear span of all functions on $A$ of the form
$$
\1_{a+F(z)},
$$
where $a\in A$ and $F$ is a total filtration on $A$, and $z\in\Z^n$ satisfies
$$
|D\cap F(z)|\ <\ \infty\qquad\mbox{and}\qquad |\hat K\cap \hat F(z)|\ <\ \infty.
$$

The sequence (\ref{S}) gives rise to a Fourier transform on $\CD_n(A)$ defined by linearity and
$$
\widehat{\1_{a+F(z)}}(\chi)\=e^{2\pi i\chi(a)}\frac{|D\cap F(z)|}{|\hat K\cap \hat F(z)|}\,\1_{\hat F(z)}(\chi).
$$
We have to check the well-definedness.
For this we suppose $z'\le z$ and that $F(z')$ has finite index in $F(z)$.
Then 
$$
\1_{F(z)}\=\sum_{a:F(z)/F(z')} \1_{a+F(z')}.
$$
We have to show that
$$
\widehat{\1_{F(z)}}\=\sum_{a:F(z)/F(z')} \widehat{\1_{a+F(z')}}.
$$
For $\chi\in\hat A$ we compute
$$
\sum_{a:F(z)/F(z')} \widehat{\1_{a+F(z')}}(\chi)\=\sum_{a:F(z)/F(z')} e^{2\pi i\chi(a)}\frac{|D\cap F(z')|}{|\hat K\cap \hat F(z')|}\,\1_{\hat F(z')}(\chi).
$$
If $\chi\notin \hat F(z)$, then the sum on the right is zero.
Otherwise it is
$$
|F(z)/F(z')|\frac{|D\cap F(z')|}{|\hat K\cap \hat F(z')|}.
$$
We have to show that this number equals $\frac{|D\cap F(z)|}{|\hat K\cap \hat F(z)|}$.
In other words, our claim is
$$
|F(z)/F(z')|\= \frac{|D\cap F(z)|}{|\hat K\cap \hat F(z)|}
\frac{|\hat K\cap \hat F(z')|}{|D\cap F(z')|}.
$$
This, however, is clear, as we have the exact sequence
$$
0\to D\cap F(z)/D\cap F(z')\to F(z)/F(z')\to F_K(z)/F_K(z')\to 0,
$$
and
\begin{eqnarray*}
|F_K(z)/F_K(z')|&=&|\widehat{F_K(Z)/F_K(z')}|\\
&=&|\hat F_K(z')/\hat F_K(z)|\\
&=& |\hat K\cap \hat F(z')/\hat K\cap \hat F(z)|.
\end{eqnarray*}
So the Fourier transform is well-defined.

\begin{theorem}
(Poisson summation formula)\\
For every $f\in\CD_n(A)$ we have
$$
\sum_{d\in D}f(d)\=\sum_{\chi\in\hat K}\hat f(\chi).
$$
Both sums are finite.
\end{theorem}

\prf
By linearity, it suffices to consider the case $f=\1_{a+F(z)}$.
In this case the left hand side of the summation formula equals
$|D\cap (a+F(z))|$.
Note that this number is zero if $a\notin D+F(z)$ and equals $|D\cap F(z)|$ otherwise.
The right hand side equals
$$
\sum_{\chi\in\hat K\cap \hat F(z)}e^{2\pi i\chi(a)} \frac{|D\cap F(z)|}{|\hat K\cap \hat F(z)|}.
$$
If $a\notin D+F(z)=(K\cap F(z))^\perp$, then this sum is zero.
Otherwise it equals $|D\cap F(z)|$.
\qed

The existence of a sequence (\ref{S}) is a non-trivial condition, as we will see in the following 2-dimensional example.
Let $R=\Z$ and let
$$
A\=\bigoplus_{r,s\in\Z} A_{r,s},
$$
where $A_{r,s}$ is a finite $\Z$-module.
The standard filtration is defined as
$$
F(j)\=\bigoplus_{\stackrel{r\le j}{s\in\Z}} A_{r,s}.
$$
For $i\le j$ there is a canonical isomorphism $F(j)/F(i)\cong \bigoplus_{\sr{i< r\le j}{s\in\Z}}A_{r,s}$.
The standard filtration on the quotient $F(j)/F(i)$ is
$$
F_{ij}(k)\=\bigoplus_{\sr{i<r\le j}{s\le k}}A_{r,s}.
$$
Then $A$ admits a sequence (\ref{S}) if and only if the following conditions are met:
\begin{enumerate}[\rm (a)]
\item $\exists j\in\Z\ \forall r>j\ \exists s_r\in\Z : s<s_r\Rightarrow A_{r,s}=0$,
\item $\exists i\in\Z\ \forall r<i\ \exists s_r\in\Z : s>s_r\Rightarrow A_{r,s}=0$.
\end{enumerate}

{\small Mathematisches Institut\\
Auf der Morgenstelle 10\\
72076 T\"ubingen\\
Germany\\
\tt deitmar@uni-tuebingen.de}

\today

\end{document}